%AMSTeX File

% masl14.tex --- 03-01-96 10:50

\documentstyle{amsppt}
\TagsOnRight
%\NoBlackBoxes
\hoffset=.1in

\define\rx{\operatorname{R}}
\define\lx{\operatorname{L}}

\define\id{\operatorname{id}}

\define\rar{ \Rightarrow }

%\define\duch{\chi_{{\phantom{}}_?}}   %?=indeks
%\define\duch{\chi_{{\phantom{}}_G}}
%\define\duch{\chi_{{\phantom{}}_Z_{nn}}}

\define\a{\alpha}
\define\be{\beta}
\define\vf{\varphi}

\define\ve{\varepsilon}

\define\de{\delta}
\define\si{\sigma}
\define\om{\omega}

\define\vd{\varDelta}
\define\vl{\varLambda}
\define\vg{\varGamma}

\define\caa{\Cal A}
\define\cab{\Cal B}

\define\cam{\Cal M}
\define\can{\Cal N}
\define\cp{\Cal P}
\define\car{\Cal R}

\define\rti{\widetilde \rho}

\define\pd#1#2{\dfrac{\partial#1}{\partial#2}}

\define\vc#1{(#1_1,\ldots,#1_n)}
\define\vct#1{[#1_1,\ldots,#1_n]}
\define\vect#1{\{#1_1,\ldots,#1_n\}}

\define\sui{\sum_{i}}
\define\suk{\sum_{k}}
\define\sul{\sum_{l}}
\define\sumka{\sum_{i,k,l}}

\define\xx{\otimes}

\define\jak{\dfrac{i}{\kappa}}

\define\rhoti{\widetilde \rho}

%\define\chi_{{\phantom{}}_{ \chi_{{\phantom{}}_} }

\define\rol{\rho_{\lx}}
\define\ror{\rho_{\rx}}
\define\tirol{\rti_{\lx}}
\define\tiror{\rti_{\rx}}

\define\lama{\vl^\mu{}_\nu}
\define\kapa{\cp_\kappa}
\define\mina{\cam_\kappa}

\define\1{$\kappa$-Poincar\'e group}
\define\2{covariant}
\define\3{calculus}
\define\4{differential}
\define\5{dimensional}
\define\6{quantum}
\define\7{quantum group}
\define\8{condition}
\define\9{generalization}
\define\0{algebra}

\define\defx{deformation}
\define\pro{properties}
\define\wor{Wo\-ro\-no\-wicz}
\define\wrt{with respect to}
\define\fol{following}
\define\bim{bimodule}

\define\kapo{$\kappa$-Poincar\'e}
\define\mink{$\kappa$-Min\-ko\-wski}
\define\linv{left-covariant}

\define\rep{representation}

\topmatter
\title Differential calculi\\
on quantum Minkowski space
\endtitle
\rightheadtext{Differential calculi}
\author Cezary Gonera$^*$, Piotr Kosi\'nski$^*$\\
{\it Department of Theoretical Physics}\\
{\it University of \L \'od\'z}\\
{\it ul. Pomorska 149/153, 90--236 \L \'od\'z, Poland}\\
Pawe\l \/ Ma\'slanka$^*$\\
{\it Department of Functional Analysis}\\
{\it University of \L \'od\'z}\\
{\it ul. St. Banacha 22, 90--238 \L \'od\'z, Poland}
\endauthor
\leftheadtext{C. Gonera, P. Kosi\'nski, P. Ma\'slanka}
\thanks
*\ \ \ Supported by KBN grant  2 P 302\,217\,06\,p\,02\newline
\endthanks

\abstract The  differential calculus on $n$-dimensional \6  Min\-ko\-wski 
space covariant with respect to left action of $\kappa$-Poincar\'e group is 
constructed  and its uniqueness is shown.
\endabstract

\endtopmatter

\document

\head I. Introduction
\endhead

The \kapo{}  algebra, introduced in [1], provides a Hopf algebra \defx{} of 
standard Poincar\'e algebra which depends on dimensionful parameter  $\kappa$. Its global 
counterpart, $\kappa$-Poincar\'e group $P_\kappa$ has been constructed 
by Za\-krze\-wski [2]. It is a free $\ast$-\0 generated by the hermitean 
elements $\lama$, $a^\mu$, subject to the \fol{} \8s
$$
\aligned
&[a^\mu,a^\nu]  = \dfrac{i}{\kappa}( \de_0^\mu a^\nu -  \de_0^\nu 
 a^\mu),\\
&[\lama,  \vl^\a{}_\be]  = 0,\\
& [\lama,  a^\a]  = -  \dfrac{i}{\kappa} ((\vl^\mu{}_0 - \de_0^\mu) 
\vl^\a{}_\nu + (\vl^0{}_\nu - \de_\nu^0)  g^{\mu\a}),\\
&\vd(\vl^\mu{}_\nu)  =  \vl^\mu{}_\a \otimes \vl^\a{}_\nu,\\
&\vd(a^\mu)  = \vl^\mu{}_\nu \otimes a^\nu + a^\mu \otimes I,\\
& S(\vl^\mu{}_\nu) =  \vl_\nu{}^\mu,\\
& S(a^\mu) = - \vl_\nu{}^\mu a^\nu,\\ 
& \ve(\vl^\mu{}_\nu)  =  \de^\mu_\nu,\\
&\ve(a^\mu) = 0.
\endaligned
\tag{1}
$$
It appears ([2], [3])  that one can also define a  noncommutative \9 of 
Min\-ko\-wski spacetime --- the \mink{} space  $\mina$. The \1 acts on 
$\mina$ \2ly from the left. Once one accepts the idea that the 
$\kappa$-deformed Poincar\'e symmetry can have something to do with reality 
the next step is to find the natural \9s of standard geometric notions 
related to Min\-ko\-wski space. The first step toward this directions was 
made by Sitarz [4] who showed that one cannot construct fourdimensional \4 
\3 on $\mina$ which is \2 \wrt{} infinitesimal left action of $\kapa$. He 
sketched also the construction of five\5 \2 \3.

In the present paper we consider the problem of the classification of \4 
calculi on $\mina$ which are \2 \wrt{} the left {\bf global} action of \1 on 
$\mina$. No restriction is  made concerning the \5ity of spacetime (i.e. the 
indices $\mu$, $\nu$, etc. in (1) run from 0 to  $n-1$). We show that the 
lowest \5 nontrivial \linv{} \3 is $n+1$-\5 and is unique. Its construction 
is given explicitly and the result coincides with the suggestion of Sitarz.

As a main tool we use the beautiful \wor{} theory of \4 calculi ([5]). In 
Section II we show that the \wor{} theory can be immediately extended to 
deal with the problem of \2 \4 calculi on \6 spaces. This natural extension 
provides a nice framework to discuss our  problem. In Section III we use the 
scheme developed in Section II to show that there is a unique $n+1$-\5 \3 on 
$\mina$ which is \linv{} \wrt{} the action of $\kapa$. This \3 is explicitly 
constructed and shown to be lowest-\5 nontrivial \3 on $\mina$.

The results obtained here were briefly reported in [6]. Let us also mention 
that the \4 calculi on \6 spacetimes \2 \wrt{} other \defx{s} of Poincar\'e 
group were considered by Podle\'s ([7]).

\head II. Covariant differential calculi
\endhead

Let us first indicate how one can extend the Woronowicz theory of \4 calculi 
([5]) to the \fol{} situation: assume that the \7 $\cab$ acts on \6 space 
$\caa$; one looks for \4 calculi on $\caa$ on which the \2 action of $\cab$ 
can be defined as 'induced' by the action of  $\cab$ on  $\caa$. All proofs 
are omitted as being a straightforward extension of those given by \wor{}.

Let $\caa$ be an \0 with unity (\6 space). The starting point in \wor{} 
construction is the universal \bim{} $\caa^2 \subset \caa \xx \caa$ defined 
by  
$$
\aligned
& \caa^2 = \{ \suk a_k \xx b_k \in \caa \xx \caa \vert  \suk a_k  b_k = 
0\},\\
& c (\suk a_k \xx b_k) =  \suk ca_k \xx b_k,\\
& (\suk a_k \xx b_k)c = \suk a_k \xx b_kc.
\endaligned
\tag{2}
$$
The universal \4 $D : \caa \to \caa^2$ is given by $da = I \xx a - a \xx I$.

It can be easily shown that any other \3 is obtained from the universal one 
by  dividing by an appropriately chosen sub-\bim{} $\can \subset \caa^2$.

Let now $\rol $ be a left action of a \7 $\cab $ on $\caa$, i.e. a 
homomorphism $\rol : \caa \to \cab \xx \caa$ obeying
$$
\aligned
& (\id \otimes \rol) \circ \rol = (\vd \xx \id) \circ  \rol,\\
& (\ve \otimes \id) \circ \rol =  \id.
\endaligned
\tag{3}
$$
Let us define the action $\tirol$ of $\cab$ on $\caa^2$ as follows
$$
\aligned
& q = \sui x_i \xx y_i \in \caa \xx \caa,\\
& \rol(x_i) =  \sum_k a_i{}^k \otimes x_i{}^k \in \cab \otimes \caa,\\
& \rol(y_i) =  \sum_l b_i{}^l \otimes y_i{}^l \in \cab \otimes \caa,
\endaligned
\tag"{(4a)}"
$$
then
$$
\tirol(q) = \sumka a_i{}^k b_i{}^l \xx x_i{}^k \xx y_i{}^l .
\tag"{(4b)}"
$$
Obviously, $\tirol : \caa \xx \caa \to \cab \xx \caa \xx \caa$, however, it 
is straightforward to show that $\tirol : \caa^2  \to  \cab \xx \caa^2$. 
Following the same lines as in [5], one easily proves the \fol{} \pro{} of 
$\tirol$
\roster
\item"{(i$_L$)}" for $x \in \caa$, \ $y \in \caa^2$
$$
\align
& \tirol (xy) = \rol(x) \tirol (y),\\
& \tirol (yx) = \tirol(y) \rol(x)
\endalign
$$
\item"{(ii$_L$)}"  
$$
\tirol \circ D = (\id \xx D) \circ \rol
\tag{5}
$$
\item"{(iii$_L$)}" 
$$
\align
& (\id \otimes \tirol) \circ \tirol = (\vd \xx \id) \circ  \tirol,\\
& (\ve \otimes \id) \circ \tirol =  \id.
\endalign
$$
\endroster

Property (iii$_L$) means that $\tirol$ is the left action of $\cab$ on 
$\caa^2$ while (i), (ii) can be summarized by saying that $\tirol$ is the 
lift of $\rol$ to $\caa^2$  ($\tirol$ is the left action of $\cab$ on 
universal \4 \3 on $\caa$ induced from the left action $\rol$). Moreover, 
let us note that the \fol{} formula holds
$$
\tirol (\sui x_i D y_i) = \sui \rol(x_i)(\id \xx D)  \rol(y_i)
\tag{6}
$$
which is a counterpart of (1.15) of [5].

Now, assume that $\can \subset \caa^2$ is a sub\bim{} such that
$$
\tirol (\can) \subset \cab \xx \can.
\tag{7}
$$

Then the \4 calculus $(\vg,d)$ determined by $\can$ has the \fol{} property
$$
\sui x_i d y_i  = 0 \rar  \sui \rol(x_i)(\id \xx d)  \rol(y_i) = 0.
\tag{8}
$$
Therefore
$$
\tirol(\sui x_i d y_i)  =  \sui \rol(x_i)(\id \xx d)  \rol(y_i)
\tag{9}
$$
is well defined linear mapping from $\vg$ into $\cab \xx \vg$. Formulae 
(i$_L$)--(iii$_L$) and (6) hold upon replacing $\caa^2$ by $\vg$ and $D$ by 
$d$.

We shall say that $(\vg,d)$ is left-\2 \wrt{} the action of $\cab$.

All the above results can be extended mutatis mutandis to right actions. Let 
$\ror : \caa \to \caa \xx \cab$ be right action of $\cab $ on $\caa$
$$
\aligned
& (\ror \otimes \id) \circ \ror =  (\id \xx \vd) \circ \ror,\\
& ( \id \otimes \ve) \circ \ror =  \id.
\endaligned
\tag{10}
$$
For
$$
\aligned
& q = \sui x_i \xx  y_i \in \caa \xx \caa,\\
& \ror (x_i) = \suk x_i{}^k \xx a_i{}^k \in \caa \xx \cab,\\
& \ror (y_i) = \sul y_i{}^l \xx b_i{}^l \in \caa \xx \cab,
\endaligned
\tag"{(11a)}"
$$
we put
$$
 \tiror (q) = \sumka x_i{}^k \xx y_i{}^l  \xx a_i{}^k  b_i{}^l .
\tag"{(11b)}"
$$
Again $\tiror : \caa^2 \to \caa^2   \xx \cab$ obeys
\roster
\item"{(i$_R$)}" for $x \in \caa$, \ $y \in \caa^2$
$$
\align
& \tiror (xy) = \ror(x) \tiror (y),\\
& \tiror (yx) = \tiror(y) \ror(x)
\endalign
$$
\item"{(ii$_R$)}"  
$$
\tiror \circ D = (D \xx \id) \circ \ror
\tag{12}
$$
\item"{(iii$_R$)}"
$$
\align
& (\tiror  \otimes \id ) \circ \tiror = (\id \xx \vd) \circ  \tiror,\\
& (\id \otimes \ve) \circ \tiror =  \id.
\endalign
$$
\endroster
as well as
$$
\tiror (\sui x_i D y_i) = \sui \ror(x_i)(D \xx \id) \ror(y_i).
\tag{13}
$$

Now assume  $\can \in \caa^2$ to be a sub-bimodule such that $\tiror(\can) 
\subset \can \otimes \cab$. Then, for the \3 $(\vg, d)$ determined by 
$\can$, (i$_R$)--(iii$_R$)  and (13) hold with appropriate replacements 
$\caa^2 \to \vg$, $D \to d$.

Finally, let the pair $(\rol,\ror)$ of actions of $\cab$ on $\caa$ be given. 
We assume that $\rol$, $\ror$ commute
$$
(\id \xx \ror) \circ \rol =  ( \rol \xx \id ) \circ \ror.
\tag{14}
$$

We say that $(\vg, d)$ is bi\2 \wrt{} the action of $\cab$ on $\caa$ if it 
is left- and right-\2. Then, it has all \pro{} of  left- and right-\2 
calculi together with the \fol{} one (cf. (1.20) of [5])
$$
(\id \xx \tiror) \circ \tirol =  ( \tirol \xx \id ) \circ \tiror.
\tag{15}
$$

Let us now discuss the problem of infinitesimal action of $\cab$ on $\caa$. 
Let $\chi$ be any element of the Hopf \0 dual to $\cab$. We put
$$
\aligned
& \chi_{{\phantom{}}_{\rol}} = (\chi \xx \id) \circ \rol,\\
& \chi_{{\phantom{}}_{\tirol}} = (\chi \xx \id) \circ \tirol.
\endaligned
\tag{16}
$$

The first definition, introduced by \wor{} ([8]), coincides with the one 
used by Majid and Ruegg ([3]). The second one is equivalent to the proposal 
of Sitarz ([4])
$$
\aligned
\chi_{{\phantom{}}_{\tirol}} (xdy) & = (\chi \xx \id) (\rol(x) (\id \xx d) 
\rol(y))\\
& = [(\chi_{{\phantom{}}_{(1)}} \xx \id) \rol(x)] (\id \xx d) 
[(\chi_{{\phantom{}}_{(2)}} \xx \id) \rol(y)]\\
& = \chi_{{\phantom{}}_{(1)\rol}}   (x) (\id \xx d)
\chi_{{\phantom{}}_{(2)\rol}} (y)
\endaligned
\tag{17}
$$
where $\vd\chi = \chi_{{\phantom{}}_{(1)}} \xx \chi_{{\phantom{}}_{(2)}}$. 
Analogous definitions can be  given for $\ror$ and $\tiror$.

From the above discussion it follows then that in order to check whether a 
\3 on $\caa$ is consistent with the action of $\cab$ on $\caa$ it is 
sufficient to check the property $\tirol (\can) \subset \cab \xx \can$ 
($\tiror (\can) \subset \can  \xx \cab$). 

This simplifies considerably if $\caa$ itself is a \6 group and $\can$ 
defines (say) \linv{} \3 on it. Then $\can = r^{-1} (\caa \xx \car)$ where 
$\car$ is a right ideal in $\ker \ve$; any element of  $\caa \xx \car$ can 
be written as
$$
t = \sui a_i \xx x_ib_i
\tag{18}
$$
where $a_i,b_i \in \caa$ and $x_i$ are generators of $\car$. From the very 
definition of the operation $r^{-1} $ the \fol{} formula follwos immediately 
([5])
$$
r^{-1} (t) = \sum_{i,l} a_i S(b^{l'}{}_i) r^{-1} (I \xx x_i) b^{l''}{}_i
\tag"{(19a)}"
$$
where
$$
\vd(b_i) = \sul b^{l'}{}_i \xx  b^{l''}{}_i.
\tag"{(19b)}"
$$

The \pro{} (5) applied to the universal \3 imply
$$
\tirol (r^{-1} (t)) = \sum_{i,l} \rol (a_i S(b^{l'}{}_i))\tirol ( r^{-1} 
(I \xx x_i) ) \rol (b^{l''}{}_i).
\tag{20}
$$
Therefore it is sufficient to check that
$$
\tirol ( r^{-1} (I \xx x_i) ) \subset \cab \xx \can
\tag{21}
$$
for all generators $x_i$ of $\car$.

Let us now pass to the external \0. Given  $\tirol : \vg \to \cab \xx \vg$ 
we define $\rhoti^{\xx 2}_{\text{L}} : \vg^{\xx 2} \to \cab \xx \vg^{\xx 2}$ 
by
$$
\rhoti^{\xx 2}_{\text{L}}(\om_1 \xx \om_2) = \sum_{k,l} a_{1k} a_{2l} \xx 
\om_{1k} \xx \om_{2l}
\tag{22}
$$
extended by linearity; here $\om_i \in \vg$ and 
$$
\tirol(\om_i) = \suk  a_{ik} \xx \om_{ik}, \qquad i = 1,2.
\tag{23}
$$

Let us assume that $\caa$ is a \7, $(\vg,d)$ --- a bi\2 \3 on it and let 
$\si$ be the module homomorphism defined in Proposition 3.1 of [5]. Then 
$\vg^{\wedge 2}$ is defined as
$$
\vg^{\wedge 2} = \vg^{\xx 2}\slash \ker(I - \si)
\tag{24}
$$
and, in order to have a consistent action of $\cab$ on $\vg^{\wedge 2}$ we 
must only check that
$$
(\id \xx \si) \circ \rhoti^{\xx 2}_{\text{L}} = \rhoti^{\xx 
2}_{\text{L}}\circ \si.
\tag{25}
$$

Due to the property
$$
\rhoti^{\xx 2}_{\text{L}} (xy) = \rol (x) \rhoti^{\xx 2}_{\text{L}}  (y), 
\qquad x \in \caa, \ \ y \in \vg^{\xx 2}
\tag{26}
$$
it is sufficient to verify (25) for the basic elements only.

\head III. Bicovariant calculi on $\kappa$-Minkowski space
\endhead

The $n$-\5  \mink{} space $\mina$  is an $\ast$-\0 with unity generated by $n$ 
hermitian elements $x^\mu$ subject to the \fol{} \8s ([2], [3])
$$
[x^\mu, x^\nu] = \jak(\de_0^\mu x^\nu   - \de_0^\nu x^\mu).
\tag{27}
$$
 $\cam_\kappa$ can be equipped with the structure of the  quantum group by 
 putting
$$
\aligned
& \vd  x^\mu = I \xx x^\mu +  x^\mu \otimes I,\\
& S(x^\mu) = - x^\mu,\\
& \ve(x^\mu) = 0.
\endaligned
\tag{28}
$$

The left action of $n$-\5 \1 $\kapa$ on $\mina$ can be defined as follows
$$
\aligned
& \rol(I) = I \otimes I,\\
& \rol(x^\mu) = \vl^\mu{}_\nu \otimes x^\nu + a^\mu \otimes I
\endaligned
\tag{29}
$$
extended by linearity and multiplicativity.

We want to find a \linv{} (\wrt{} action of $\kapa$) calculi on $\mina$. As 
the first step let us note that $\mina$ is a subgroup of $\kapa$. Indeed, 
$\Pi : \kapa \to \mina$ given by
$$
\Pi(a^\mu) = x^\mu, \qquad \Pi(\lama) = \de^\mu_\nu I
\tag{30}
$$
is an  epimorphism obeying
$$
\vd_\cam \circ \Pi = (\Pi \xx \Pi) \circ \vd_\cp.
\tag"{(31a)}"
$$
Moreover, it is immediate to check that
$$
(\Pi \otimes \id)   \circ \rol = \vd_\cam.
\tag"{(31b)}"
$$

Let $\tirol$ be the  extention of $\rol$ to $\cam_\kappa^2$. Equations (8), 
(31) and the results contained in [5] imply that any \3 on $\mina$ \linv{} 
\wrt{} action of $\kapa$ is also \linv{} \wrt{} action of $\mina$ on itself. 
Therefore the relevant sub-\bim{} $\can$ is of the form $r^{-1}(\mina \xx 
\car)$ where $\car$ is a right ideal in $\ker \ve_\cam$.

Let $\car$ be any ideal in $\ker \ve_\cam$. Any $a \in \car$ can be written 
as
$$
a = \sum_{\mu_0,\mu_k} c_{\underline{\mu}} (x^0)^{\mu_0} \prod^{n-1}_{k=1} (x^k)^{\mu_k}.
\tag{32}
$$

Let us call $|\mu| = \mu_0 + \sum^{n-1}_{k=1} \mu_k$; obviously, 
$c_{\underline{\mu}} = 0$  for  $|\mu| = 0$; further, let
$$
\aligned
& \mu(a) = \max_{c_{\underline{\mu}} \ne 0} |\mu|,\\
& \mu(\car)  = \min_{a \in \car} \mu(a).
\endaligned
\tag{33}
$$

Obviously, $\mu(\car) \ge 1$; let us first assume that $\mu(\car) = 1$. This 
means that $c_0x^0 + c_kx^k \in \car$ for some (not all zero) constants 
$c_0$, $c_k$. But
$$
\tirol (r^{-1}(I \xx c_\mu x^\mu)) = c_\mu \lama \xx  r^{-1}(I \xx x^\nu).
\tag{34}
$$
Therefore $x^\mu \in \car$ for all $\mu$, i.e. $\car = \ker \ve_\cam$ and 
the corresponding \3 is trivial.

As the next step let us take $\mu(\car) = 2$. It is straightforward to check 
that
$$
\tirol (r^{-1}(I \xx  x^{\mu\nu})) = \vl^\mu{}_\a \vl^\nu{}_\be \xx  
r^{-1}(I \xx x^{\mu\nu})
\tag{35}
$$
where
$$
 x^{\mu\nu} \equiv  x^{\mu}x^{\nu} + \jak (g^{\mu\nu} x^0 -  g^{0\mu}x^{\nu}).
\tag{36}
$$

Due to the fact that $\vl$'s commute among themselves we can write a 
standard \rep{} theory of Lorentz group. First of all, we note that 
$r^{-1}(I \xx  x^{\mu\nu})$ transform as a second order symmetric 
($x^{\mu\nu} = x^{\nu\mu}$ due to (27), (36)) tensor. It carries $D^{(1,1)} 
\oplus D^{(0,0)} $ \rep{} of Lorentz group. Let us first take all 
$x^{\mu\nu}$ as generators of $\car$. Then $(x^0)^2 = x^{00} \in \car$ and 
$x^ix^0 = x^{i0} \in \car$, i.e. $[(x^0)^2 , x^i] \in \car$; therefore 
$x^i x^0 + x^0 x^i \in \car$ and $x^0 x^i \in \car$. However, $x^{0i} 
= x^0 x^i - \jak x^i \in \car$ which implies  $x^i \in \car$. Then 
Poincar\'e invariance implies $x^0 \in \car$ and  our \3 is trivial.

To improve the situation we can only, due to \8 (21), subtract $D^{(0,0)}$ 
or  $D^{(1,1)} $. Obviously, subtracting  $D^{(1,1)} $ gives larger \3, so we 
will substract $D^{(0,0)} $. It is not difficult to check then that the 
\fol{} lemma holds.

\proclaim{Lemma} Let $\car \subset \ker \ve_\cam$ be right ideal generated 
by the elements
$$
x^{\mu}x^{\nu} + \jak(g^{\mu\nu} x^0 -  g^{0\mu}x^{\nu}) - \dfrac{1}{n} 
g^{\mu\nu} \Big(x^2 + \dfrac{i(n-1)}{\kappa} x^0\Big).
\tag{37}
$$
Then
\roster
\item"{(a)}" $\car$ defines a left-$\kapa$-\2 \3 on $\mina$,
\item"{(b)}" $a \in \car$ implies $S(a)^* \in \car$,
\item"{(c)}" $ \ker \ve_\cam\slash \car$  is spanned by $x^\mu$ and by
$$
\vf \equiv x^2 + \jak (n-1) x^0.
\tag{38}
$$
\endroster
\endproclaim

Now, we can construct the relevant \3. The  left-invariant forms are 
$$
\aligned
& \tau^\mu = \pi r^{-1} (I \xx x^\mu) = dx^\mu,\\
& \tau = \pi r^{-1}(I \otimes \vf) = d\vf - 2x_\mu dx^\mu
\endaligned
\tag{39}
$$
and they appear to be also right invariant. The commutation rules are easily 
derived according to the standard procedure of [5]
$$
\aligned
& [\tau^\mu, x^\nu] = \jak g^{0\mu} \tau^{\nu} - \jak g^{\mu\nu} \tau^0 + 
\dfrac{1}{n}  g^{\mu\nu} \tau,\\
& [\tau,x^\mu]  = - \dfrac{n}{\kappa^2}  \tau^{\mu} 
\endaligned
\tag{40}
$$
while the hermicity \pro{} read
$$
\aligned
& (\tau^\mu)^*  =  \tau^\mu,\\
& \tau^*  = -   \tau.
\endaligned
\tag{41}
$$
The left action of $\kapa$ on $\mina$ is easily calculated to be
$$
\aligned
& \tirol (\tau^\mu)  = \lama \xx \tau^\nu,\\
& \tirol  (\tau)  = I \xx \tau.
\endaligned
\tag{42}
$$

In order to construct the external \0 we first  verify property (25) for the 
\bim{} homomorphism $\si$: \ $\si(\tau^\mu \xx \tau^\nu) = \tau^\nu \xx 
\tau^\mu$, $\si(\tau \xx \tau^\mu) =\tau^\mu \xx \tau$, $\si(\tau^\mu \xx 
\tau) = \tau \xx \tau^\mu$. The external \0 implied by $\si$ takes the 
standard form
$$
\aligned
& \tau^\mu \wedge  \tau^\nu =  -\tau^\nu \wedge \tau^\mu     ,\\
& \tau \wedge \tau^\mu   = -  \tau^\mu \wedge  \tau.
\endaligned
\tag{43}
$$
Moreover,
$$
\aligned
& d\tau^\mu = 0   ,\\
& d\tau = -2d \tau^\mu  \wedge d \tau_\mu.
\endaligned
\tag{44}
$$

From the discussion carried out above it follows that the $n+1$-\5 \3 
described by (39)--(43) is the lowest \5 nontrivial \3 on $\mina$ \2 \wrt{} 
the left action of $\kapa$. This is due to the fact that all \4 calculi with 
$\mu(\car) \ge 3$ have higher dimensions.

Finally, let us compare our \3 with that proposed by Sitarz ([4]). Our 
equations (40) agree with equations (60) of [4] under the identification: \ 
$x^\mu \to ix^\mu$, $\tau \to \frac{4}{\kappa^2}\vf$. In the two\5 case 
there is also in agreement provided the replacement $\tau \to 
\frac{2}{\kappa^2}\vf$ is made; also the multiplication rules for one-forms 
(equations (58) of [4]) coincide in this case.
\vskip.2cm
{\bf  ACKNOWLEDGMENTS.}
\flushpar
We gratefully acknowledge the collaboration of Jan Sobczyk at the early 
stages of this work.

\enddocument